\documentclass{article}
\usepackage{amsmath,amssymb,amsthm,amscd,dsfont}
\numberwithin{equation}{section} \allowdisplaybreaks
\begin{document}
\newtheorem{theorem}{Theorem}[section]
\newtheorem{defin}{Definition}[section]
\newtheorem{prop}{Proposition}[section]
\newtheorem{corol}{Corollary}[section]
\newtheorem{lemma}{Lemma}[section]
\newtheorem{rem}{Remark}[section]
\newtheorem{example}{Example}[section]
\title{Geometric quantization of weak-Hamiltonian functions}
\author{{\small by}\vspace{2mm}\\Izu Vaisman}
\date{}
\maketitle
{\def\thefootnote{*}\footnotetext[1]%
{{\it 2000 Mathematics Subject Classification:} 53D17, 53D50.
\newline\indent{\it Key words and phrases}: big-isotropic structure;
geometric prequantization data system; prequantization space;
polarization.}}
\begin{center} \begin{minipage}{12cm}
A{\footnotesize BSTRACT. The paper presents an extension of the
geometric quantization procedure to integrable, big-isotropic
structures. We obtain a generalization of the cohomology
integrality condition, we discuss geometric structures on the
total space of the corresponding principal circle bundle and we
extend the notion of a polarization.}
\end{minipage} \end{center}
\vspace*{5mm}
\section{Big-isotropic structures}
Weak-Hamiltonian functions belong to the framework of big-isotropic
structures and have been discussed in \cite{{V-iso},{V-JMP}}. For
the convenience of the reader, we recall some basic facts here.

All the manifolds and mappings are of $C^\infty$ class and we denote
by $M$ an $m$-dimensional manifold, by $\chi^k(M)$ the space of
$k$-vector fields, by $\Omega^k(M)$ the space of differential
$k$-forms, by $\Gamma$ the space of global cross sections of a
vector bundle, by $X,Y,..$ either contravariant vectors or vector
fields, by $\alpha,\beta,...$ either covariant vectors or $1$-forms,
by $d$ the exterior differential and by $L$ the Lie derivative.

The vector bundle $T^{big}M=TM\oplus T^*M$ is called the {\it big
tangent bundle}. It has the natural, non degenerate metric of zero
signature (neutral metric)
\begin{equation}\label{gFinC}
g((X,\alpha),(Y,\beta))=\frac{1}{2}(\alpha(Y)+\beta(X)),
\end{equation} the non degenerate, skew-symmetric $2$-form
\begin{equation}\label{omegainC}
\omega((X,\alpha),(Y,\beta))= \frac{1}{2}(\alpha(Y)-\beta(X))
\end{equation} and the Courant bracket of cross sections
\begin{equation}\label{crosetC} [(X,\alpha),(Y,\beta)]_C = ([X,Y],
L_X\beta-L_Y\alpha+\frac{1}{2}d(\alpha(Y)-\beta(X)));\end{equation}
unless to avoid confusion, the index $C$ will be omitted.
\begin{defin}\label{defbigiso} {\rm A $g$-isotropic subbundle
$E\subseteq T^{big}M$ of rank $k$ $(0\leq k\leq m)$ is called a {\it
big-isotropic structure} on $M$. A big-isotropic structure $E$ is
{\it integrable} if $\Gamma E$ is closed by Courant brackets.}
\end{defin}

From the properties of the Courant bracket it follows that the
integrability of $E$ is equivalent with the property: $$[\Gamma
E,\Gamma E']_C\subseteq\Gamma E',$$ where $E'\perp_g E$,
$(E\subseteq E')$ \cite{{V-iso},{V-JMP}}. The big-isotropic
structures are a generalization of the almost Dirac structures
(Dirac structures, in the integrable case), which are obtained for
$k=m$.
\begin{example}\label{exgraphlambda} {\rm\cite{V-iso} Let $S\subseteq TM$ be
a $k$-dimensional, regular distribution on $M$ and
$\lambda\in\Omega^2(M)$. Put
\begin{equation}\label{Elambda}
E_{(\lambda,S)}=graph(\flat_\lambda|_{S})=\{(X,\flat_{\lambda}X
=i(X)\lambda)\,/\,X\in
S\}.\end{equation} Then $E_{(\lambda,S)}$ is a big isotropic
structure with the $g$-orthogonal bundle
\begin{equation}\label{perpElambda} E'_{(\lambda,S)}=\{(Y,\flat_{\lambda}
Y+\gamma)\,/\,Y\in TM,\,\gamma\in{\rm ann}\,S\}.\end{equation} If
$d\lambda=0$, $E_{(\lambda,S)}$ is integrable iff $S$ is a
foliation.}\end{example}
\begin{example}\label{exgraphP} {\rm\cite{V-iso} Let $\Sigma$
be a rank $k$ subbundle of $T^*M$ and $P\in\chi^2(M)$. Then
\begin{equation}\label{eqEP} E_{(P,\Sigma)}=graph(\sharp_P|_{\Sigma}) =
\{(\sharp_P\sigma=i(\sigma)P,\sigma)\,/\,\sigma\in
\Sigma\}\end{equation} is a big-isotropic structure on $M$ with the
$g$-orthogonal bundle
\begin{equation}\label{eqE'P} E'_{(P,\Sigma)}
=\{(\sharp_P\beta+Y,\beta)\,/\,\beta\in T^*M, Y\in {\rm
ann}\,\Sigma\}.\end{equation} If $P$ is a Poisson bivector field the
structure (\ref{eqEP}) is integrable iff $\Sigma$ is closed with
respect to the bracket of $1$-forms defined by
\begin{equation}\label{croset1forme} \{\alpha,\beta\}_P=
L_{\sharp_P\alpha}\beta-L_{\sharp_P\beta}\alpha
-d(P(\alpha,\beta)).\end{equation}}\end{example}

For geometric quantization it is important to point out the
existence of an adequate cohomology associated with an integrable
big-isotropic structure $E$ \cite{V-iso}. In the formulas below and
in the remaining part of the paper, calligraphic letters denote
pairs, $\mathcal{X}=(X,\alpha)$, $ \mathcal{Y}=(Y,\beta)$, etc. The
cochain spaces are the spaces of {\it truncated forms}
\begin{equation}\label{cochainEE'}
\Omega_{tr}^s(E)=\{\lambda:\wedge^{s-1}E\otimes E'\rightarrow
C^\infty(M)\,/\,\lambda|_{\wedge^{s-1}E\otimes E}\in\wedge^s
E^*\}.\end{equation} The coboundary operator is defined by
$$d_{tr}\lambda( \mathcal{X}_1,...,\mathcal{X}_s,\mathcal{Y})=
\sum_{a=1}^s(-1)^{a+1}X_a(\lambda(
\mathcal{X}_1,...,\hat{\mathcal{X}}_a,...,\mathcal{X}_s,\mathcal{Y}))$$
$$+(-1)^sY(\lambda( \mathcal{X}_1,...,\mathcal{X}_s))
+\sum_{a<b=1}^s(-1)^{a+b}\lambda( [\mathcal{X}_a,\mathcal{X}_b]_C,
\mathcal{X}_1,...,\hat{\mathcal{X}}_a,...,\hat{\mathcal{X}}_b,...,
\mathcal{X}_s,\mathcal{Y})$$ $$-\sum_{a=1}^s(-1)^a\lambda(
\mathcal{X}_1,...,\hat{\mathcal{X}}_a,...,\mathcal{X}_s,
[\mathcal{X}_s,\mathcal{Y}]_C),$$ where $\mathcal{X}\in\Gamma
E,\mathcal{Y}\in\Gamma E'$. Since all the arguments that we use are
pair-wise $g$-orthogonal, the Courant brackets above have the
properties of the bracket of a Lie algebroid (the Jacobi identity in
particular), therefore, the coboundary condition $d_{tr}^2=0$ holds.
The corresponding cohomology spaces, called {\it truncated
cohomology spaces}, will be denoted by $H_{tr}^s(E)$. The mapping
$j:\Omega^s(M)\rightarrow\Omega_{tr}^s(E)$ defined by
$$(j\lambda)(\mathcal{X}_1,...,\mathcal{X}_{s-1},\mathcal{Y})=
\lambda(X_1,...,X_{s-1},Y)$$ is a morphism of cochain complexes,
hence, there exist induced homomorphisms
$$j^*:H^s_{deR}(M,\mathds{R})\rightarrow H_{tr}^s(E).$$
\begin{rem}\label{cohgen} {\rm In an appendix to this paper we
will show that truncated cohomology has a further generalization
in the case of a pair of Lie algebroids.}\end{rem}

For the integrable big-isotropic structure $E$ the form
(\ref{omegainC}) defines the truncated $2$-form
$\omega_E=\omega|_{E\times E'}$ and a straightforward calculation
gives $d_{tr}\omega_E=0$, hence, one has a {\it fundamental
cohomology class} $[\omega_E]\in H_{tr}^2(E)$ \cite{V-iso}.

We end this section by explaining what is meant by a
weak-Hamiltonian function with respect to an integrable,
big-isotropic structure $E$ (see details in \cite{V-iso}).
\begin{defin}\label{defiweakH} {\rm
A function $f\in C^\infty(M)$ is a {\it Hamiltonian}, respectively
{\it weak-Hamiltonian}, function if there exists a vector field
$X_f\in\chi^1(M)$ such that $(X_f,df)\in\Gamma E$, respectively
$(X_f,df)\in\Gamma E'$. The vector field $X_f$ is a {\it
Hamiltonian}, respectively {\it weak-Hamiltonian}, vector field of
$f$.}\end{defin}

The vector field $X_f$ is required to be differentiable but it may
not be unique. In the Hamiltonian case $X_f$ is defined up to the
addition of any $Z\in\chi^1(M)\cap\Gamma E$ and in the
weak-Hamiltonian case up to $Z\in\chi^1(M)\cap\Gamma E'$. We denote
by $C^\infty_{Ham}(M,E)$ the set of Hamiltonian functions, by
$C^\infty_{wHam}(M,E)$ the set of weak-Hamiltonian functions and by
$\mathcal{X}_f=(X_f,df)$ the pairs described in Definition
\ref{defiweakH}.

Furthermore, if $f\in C^\infty_{Ham}(M,E)$ and $h\in
C^\infty_{wHam}(M,E)$ one defines the {\it Poisson bracket}
\begin{equation}\label{crosetfh}
\{f,h\}=-\omega_E(\mathcal{X}_f,\mathcal{X}_h)=X_{f}h=-X_{h}f,\end{equation}
which is easily seen not to depend on the choice of the Hamiltonian
vector fields. The integrability of $E$ implies $\{f,h\}\in
C^\infty_{wHam}(M,E)$ and shows that one of the weak-Hamiltonian
vector fields of the function $\{f,h\}$ is the Lie bracket
$[X_f,X_h]$. If both $f,h\in C^\infty_{Ham}(M,E)$, the Poisson
bracket is skew symmetric and belongs to $C^\infty_{Ham}(M,E)$. The
Poisson bracket satisfies the {\it Leibniz rule}
\begin{equation}\label{Leibniz} \{l,\{f,h\}\} =
\{\{l,f\},h\}\}+\{f,\{l,h\}\}, \end{equation} $\forall l,f\in
C^\infty_{Ham}(M,E),h\in C^\infty_{wHam}(M,E)$. Equality
(\ref{Leibniz}) restricts to the Jacobi identity on
$C^\infty_{Ham}(M,E)$.

In the case of a Dirac structure discussed in \cite{C}, since
$E'=E$, the notions of weak-Hamiltonian and Hamiltonian functions
coincide.
\begin{rem}\label{obsCRF} {\rm The notion of an integrable
big-isotropic structure may be complexified and all the previous
result hold if complex values are assumed overall. Important
examples are offered by the generalized CRF-structures, where $E$
appears as the $i$-eigenbundle of a skew-symmetric endomorphism
$\mathcal{F}$ of $T^{big}M$ such that $\mathcal{F}^3+\mathcal{F}=0$
\cite{V-CRF} and, in the Dirac case, by generalized complex
structures, i.e., integrable $i$-eigenbundles of a skew-symmetric
endomorphism $\mathcal{I}$ of $T^{big}M$ such that
$\mathcal{I}^2=-Id$ \cite{Galt}.}\end{rem}
\section{Prequantization of big-isotropic structures}
Prequantization is the first step of the geometric quantization
procedure. It was defined and studied, independently, by J. M.
Souriau \cite{Sour} and B. Kostant \cite{K} for symplectic
manifolds, extended by several authors (e.g., \cite{V-carte}) to
Poisson manifolds and by A. Weinstein and M. Zambon to Dirac
manifolds \cite{WZ}. Here, we extend prequantization further, to
integrable big-isotropic structures. The path to follow is
directly suggested by the Poisson and Dirac case.

Let $E$ be an integrable big-isotropic structure on $M$ and
consider a triple $(K,\nabla,\theta)$ where $K$ is a Hermitian
line bundle on $M$, $\nabla$ is a Hermitian connection on $K$ and
$\theta$ is a truncated $1$-cochain of $E$. Notice that the
isomorphism $\sharp_g$ (defined like in Riemannian geometry)
yields an isomorphism $\Omega_{tr}^1(E)=E^{'*}\approx T^{big}M/E$,
hence, the cochain $\theta$ may be seen as a pair
$(U,\nu)\in\Gamma T^{big}M$ defined up to the addition of any
$(X,\alpha)\in\Gamma E$. Our convention for this identification
will be
\begin{equation}\label{calculcociclu}
\theta(Y,\beta)=\nu(Y)+\beta(U)\hspace{3mm}((Y,\beta)\in
E').\end{equation}

The triple $(K,\nabla,\theta)$ will be called a {\it g.p. (geometric
prequantization) data system} if the {\it modified Kostant-Souriau
formula}
\begin{equation}\label{KS} \hat
hs=\nabla_{X_h}s+2\pi i(\theta( \mathcal{X}_h)+h)s\;\;(s\in\Gamma
K,h\in C^\infty_{wHam}(M,E))\end{equation} associates with every
weak-Hamiltonian function $h$ and every weak-Hamiltonian vector
field $X_h$ an operator $\hat h:\Gamma K\rightarrow\Gamma K$ with
the following property: $\forall f\in C^\infty_{Ham}(M,E), h\in
C^\infty_{wHam}(M,E)$ the commutant $[\hat f,\hat h]=\hat f\circ\hat
h-\hat h\circ\hat f$ is equal to the operator $\widehat{\{f,h\}}$
associated to the Poisson bracket and to the choice
$X_{\{f,h\}}=[X_f,X_h]$ of the weak-Hamiltonian vector field of
$\{f,h\}$. In physics, the terminology is: {\it observable} for the
function $h$ and {\it quantum operator} for the operator $\hat h$
and for other operators with a similar role.
\begin{prop}\label{condcurbpreq} A triple
$(K,\nabla,\theta)$ is a g.p. data system for $E$ iff the
curvature of the connection $\nabla$ satisfies the following
condition
\begin{equation}\label{curburapct} R_\nabla(X,Y)=2\pi i
(\omega_E(\mathcal{X},\mathcal{Y})-
(d_{tr}\theta)(\mathcal{X},\mathcal{Y})),\end{equation}
$\forall\mathcal{X}=(X,\alpha)\in E_x,\mathcal{Y}=(Y,\beta)\in
E'_x,x\in M$.\end{prop}
\begin{proof} A straightforward calculation gives
$$[\hat f,\hat h](s)=\widehat{\{f,h\}}s +
R_\nabla(X_f,X_h)s+2\pi
i[(d_{tr}\theta)(\mathcal{X}_f,\mathcal{X}_h) +\{f,h\}]s,$$ where
$\widehat{\{f,h\}}$ is constructed with the weak-Hamiltonian vector
field $X_{\{f,h\}}=[X_f,X_h]$. Accordingly, we have a g.p. data
system iff
\begin{equation}\label{condptcurbura}R_\nabla(X_f,X_h)=
-2\pi i(\{f,h\}+
d_{tr}\theta(\mathcal{X}_f,\mathcal{X}_h)),\end{equation} $\forall
f\in C^\infty_{Ham}(M,E), h\in C^\infty_{wHam}(M,E)$, which is
equivalent with the point-wise condition (\ref{curburapct}). The
domains of the arguments in (\ref{curburapct}) are explained by the
fact that $\forall(X,\alpha)\in E_x,(Y,\beta)\in E'_x$ there are
functions $f\in C^\infty_{Ham}(M,E),h\in C^\infty_{wHam}(M,E)$ such
that $d_xf=\alpha,d_xh=\beta,X_f(x)=X,X_h(x)=Y$.
\end{proof}

The following result is an easy consequence of (\ref{curburapct}).
\begin{prop}\label{integralitate} The integrable big-isotropic
structure $E$ admits g.p. data systems iff there exists a closed
$2$-form $\Theta$ on $M$ which represents an integral de Rham
cohomology class $[\Theta]\in H^2_{deR}(M)$ such that
$j^*[\Theta]=[\omega_E]$.\end{prop}
\begin{proof} Since the $2$-form
$-(1/2\pi i)R_\nabla$ represents the first Chern class of $K$,
(\ref{curburapct}) implies the required conclusion. Conversely, it
is well known that if $-\Theta\in\Omega^2(M)$ represents an
integral cohomology class $c\in H^2(M,\mathds{Z})$ there exists a
Hermitian line bundle $K$ with the first Chern class $c$ endowed
with a Hermitian connection $\nabla$ of curvature $2\pi i\Theta$.
Furthermore, if $j^*[\Theta]=[\omega_E]$ there exists a
$1$-cochain $\theta\in\Omega^1_{tr}(M,E)$ such that
\begin{equation}\label{auxintegral}
\Theta(X,Y)=\omega_E(\mathcal{X},\mathcal{Y})
-(d_{tr}\theta)(\mathcal{X},\mathcal{Y}),\end{equation} which
precisely is (\ref{curburapct}).\end{proof}
\begin{rem}\label{obsclasificare} {\rm
We may refer to the classification of the set of g.p. data systems
like in \cite{WZ}. Fix a cohomology class $c\in H^2(M,\mathds{Z})$
with image $-[\omega_E]\in H_{tr}^2(M,E)$; then, the g.p. data
systems $(K,\nabla,\theta)$ where the line bundle $K$ has first
Chern class $c$ are said to have topological type $c$. If the
topological type is fixed, $K$ is determined up to an isomorphism.
Moreover, we may also consider that the corresponding $2$-form
$\Theta$ $([\Theta]=-c)$, therefore the connection $\nabla$ too is
fixed up to the previous isomorphism. Indeed, a change
$\Theta\mapsto\Theta+d\xi$ is equivalent with a change of the
cochain $\theta\mapsto\theta+j(\xi)$. But, if $\Theta$ is fixed,
all the corresponding possible g.p. data systems are produced by
all the $1$-cochains $\theta$ that satisfy (\ref{auxintegral}).
Since condition (\ref{auxintegral}) implies that $d_{tr}\theta$ is
well defined, $\theta$ itself is defined up to the addition of a
$1$-cocycle $\kappa\in\Omega_{tr}^1(E)$, $d_{tr}\kappa=0$.
Therefore, the set of isomorphism classes of g.p. data systems of
the topological type $c$ is in a bijective correspondence with the
set of $d_{tr}$-closed truncated $1$-forms.}\end{rem}

The expression (\ref{auxintegral}) of the integrality condition
may be put into the following simpler form.
\begin{prop}\label{simpleintegrality} The integrable big-isotropic
structure $E$ admits g.p. data systems iff there exists a closed
$2$-form $\Xi$ on $M$ that represents an integral de Rham
cohomology class $[\Xi]\in H^2_{deR}(M)$ and a vector field
$U\in\chi^1(M)$ such that
\begin{equation}\label{integralcuU}
\beta(X)+(L_U\beta)(X)-\alpha([Y,U])=\Xi(X,Y),\hspace{2mm}
\forall(X,\alpha)\in\Gamma E,(Y,\beta)\in\Gamma E'.
\end{equation}
\end{prop}
\begin{proof} If the differential $d_{tr}\theta$
of the right hand side of (\ref{auxintegral}) is replaced by its
expression, while keeping in mind that
$\mathcal{X}\perp_g\mathcal{Y}$, then (\ref{auxintegral}) becomes
(\ref{integralcuU}) with $\Xi=-(\Theta+d\nu)$.\end{proof}
\begin{rem}\label{obsLU} {\rm If we denote
$L_U\mathcal{Y}=(L_UY,L_U\beta)$ condition (\ref{integralcuU}) may
be written under the form \begin{equation}\label{integralcuLU}
\omega_E( \mathcal{X},\mathcal{Y})+2g(
\mathcal{X},L_U\mathcal{Y})=\Xi(X,Y),\hspace{3mm}
\forall\mathcal{X}\in\Gamma E,\mathcal{Y}\in\Gamma
E'.\end{equation}}\end{rem}

Proposition \ref{simpleintegrality} shows that the $1$-form $\nu$
is not essential for prequantization. In fact, we have
\begin{prop}\label{corolarnu0} If $(K,\nabla,\theta)$, where $\theta$
is defined by (\ref{calculcociclu}), is a g.p. data system for the
integrable big-isotropic structure $E$ then $(K,\nabla',\theta')$,
where $\theta'(Y,\beta)=\beta(U)$, $\nabla'=\nabla+2\pi i\nu$,
also is a g.p. data system of $E$, which yields the same quantum
operators (\ref{KS}) under the form
\begin{equation}\label{KSredus} \hat
hs=\nabla'_{X_h}s+2\pi i(U(h)+h)s\;\;(s\in\Gamma K,h\in
C^\infty_{wHam}(M,E)).\end{equation} Moreover,
$\forall\tilde{\nu}\in\Omega^1(M)$, $(K,\tilde{\nabla}=\nabla+2\pi
i(\nu-\tilde{\nu}),
\tilde{\theta}(\mathcal{Y})=\beta(U)+\tilde{\nu}(Y))$ also is a g.p.
data system with the same quantum operators.\end{prop}
\begin{proof} The new triples satisfy condition
(\ref{curburapct}).\end{proof}

A cochain of the type $\theta'(Y,\beta)=\beta(U)$ will be called a
{\it vectorial cochain} and Proposition \ref{corolarnu0} shows
that it suffices to consider g.p. systems with vectorial cochains
only. However, we will continue to write the formulas for
arbitrary complex cochains $\theta$.

In the case of a Dirac structure the conditions stated in
Propositions \ref{condcurbpreq}, \ref{integralitate} coincide with
those given in \cite{WZ}. Below, we discuss the prequantization
condition in Examples (\ref{exgraphlambda}) and (\ref{exgraphP}).
\begin{example}\label{condqex1} {\rm Let $E_{(\lambda,S)}$ be the
integrable big-isotropic structure associated with the closed
$2$-form $\lambda$ and the foliation $S$ of $M$. Then, the
prequantization condition is (\ref{integralcuU}) where
$$X\in S,\,\alpha=\flat_\lambda X,\,Y\in TM,\,\beta=\flat_\lambda
Y+\gamma,\,
\gamma\in{\rm ann}\,S.$$ The corresponding result is
$$(\lambda+L_U\lambda)(X,Y)+\gamma([U,X])=-\Xi(X,Y).$$
The case $Y=0$ shows that, $\forall X\in S$, $[X,U]\in S$, i.e.,
$U$ must be projectable onto the space of leaves of the foliation
$S$. Furthermore, since
$L_U\lambda=i(U)d\lambda+di(U)\lambda=di(U)\lambda$, the
prequantization condition reduces to the fact that $\lambda$ is a
closed, integral $2$-form. In particular, we see that
$E_{(\lambda,S)}$ is prequantizable for any $S$ using the cochain
$\theta=0$. The classical case of a symplectic manifold is
included here.}\end{example}
\begin{example}\label{condex2} {\rm Let $E_{(P,\Sigma)}$ be the
integrable big-isotropic structure associated with the Poisson
bivector field $P$ and the $\{\,,\,\}_P$-closed subbundle
$\Sigma\subseteq T^*M$. The prequantization condition is
(\ref{integralcuU}) where
$$\alpha\in\Sigma,\,X=\sharp_P\alpha,\,
\beta\in T^*M,\,Y=\sharp_P\beta+Z,\,Z\in{\rm ann}\,\Sigma$$ and a
straightforward calculation gives
\begin{equation}\label{integralptPSigma}
P(\alpha,\beta)-(L_UP)(\alpha,\beta)+\alpha([U,Z])=\Xi(\sharp_P\alpha,
\sharp_P\beta+Z),\end{equation} where $U\in\chi^1(M)$ and $\Xi$
is a closed, integral $2$-form on $M$. For a Poisson structure
$\Sigma=T^*M$ and $Z=0$. Accordingly, (\ref{integralptPSigma})
reduces to the known prequantization condition of a Poisson
structure
\cite{V-carte}. Furthermore, if $P$ is defined by an integral
symplectic form, the integrality condition is satisfied for any
$\Sigma$ if we take $\Xi$ equal to the symplectic form and
$U=0$.}\end{example}
\begin{rem}\label{prodscalar} {\rm
The quantum operators of physics act on a Hilbert space. We
indicate the following procedure to transfer prequantization to a
pre-Hilbert space \cite{{V-Torino},{V-carte}}; then, a
corresponding Hilbert space can be constructed by completion. A
complex half-density is a geometric object $\rho$ with one complex
component $\rho_\alpha$ with respect to local coordinates
$(x^i_\alpha)$ on the coordinate neighborhood $U_\alpha$ such that
on $U_\alpha\cap U_\beta$ one has
$$\rho_\beta=|det(\partial x^i_\alpha/\partial
x^j_\beta)|^{1/2}\rho_\alpha.$$ The complex half-densities define
a line bundle $S^{1/2}$ on $M$. The Lie derivative acts on
half-densities by
$$(L_X\rho)_\alpha=\xi^i_\alpha\frac{\partial \xi^i_\alpha}{\partial
x^i_\alpha}+\frac{1}{2}\rho_\alpha\frac{\partial
\rho_\alpha}{\partial
x^i_\alpha}\hspace{5mm}(X=\xi^i\frac{\partial}{\partial x^i}),$$
where the Einstein summation convention is used. The quantum
operators defined by (\ref{KS}) extend to $\Gamma(K\otimes
S^{1/2}(M))$ by putting
\begin{equation}\label{hathat}\hat{\hat h}(s\otimes\rho)=(\hat
hs)\otimes\rho+s\otimes(L_{X_h}\rho).\end{equation} Using bases of
$K$ and $S^{1/2}$, we see that any cross section
$\sigma\in\Gamma(K\otimes S^{1/2}(M))$ has representations
$\sigma=s\otimes\rho$ and that $\hat{\hat h}\sigma$ is independent
of the choice of the representation. Then, the space
$^c\Gamma(K\otimes S^{1/2}(M))$, where the index $c$ means that we
take cross sections with a compact support, has the natural scalar
product
\begin{equation}\label{prodcudens}<s_1\otimes\rho_1,s_2\otimes\rho_2>
=\int_M<s_1,s_2>\rho_1\bar{\rho_2},\end{equation} where
$<s_1,s_2>$ is the Hermitian scalar product on $K$. The density
version of Stoke's theorem (e.g., \cite{V-Torino}) leads to the
fact that the operators $\hat{\hat h}$ are antiunitary with
respect to the product (\ref{prodcudens}).}\end{rem}
\begin{rem}\label{obspreqcomlex} {\rm The results on prequantization also hold
for complex big-isotropic structures, with the difference that the
cohomology classes $[\omega_E],[\Theta],[\Xi]$ of the integrality
condition are complex cohomology classes.}\end{rem}
\section{The prequantization space}
Let $E$ be an integrable big-isotropic structure on $M$ and
$(K,\nabla,\theta)$ a g.p. data system. Let $p:Q\rightarrow M$ be
a principal circle bundle such that $K$ is associated to $Q$.
Following \cite{WZ}, the total space $Q$ will be called the {\it
prequantization space}. In the case of a symplectic manifold, the
prequantization space is a contact manifold and it was the basic
object in Souriau's version of geometric quantization \cite{Sour}.
If $\|\,\|$ is the Hermitian norm on $K$, we may take
\begin{equation}\label{spQ} Q=\{b\in K\,/\,\|b\|=1\}.\end{equation}

It is known that one has the following important geometric
elements: 1) the vertical vector field $V\in\chi^1(Q)$ defined by
the infinitesimal action of a basis of the Lie algebra $
\mathit{u}(1)$ of $S^1$ by right translations, 2) the $1$-form
$\sigma\in\Omega^1(Q)$ of the principal bundle connection on $Q$
that is equivalent with the covariant derivative $\nabla$; the
form $\sigma$ vanishes on vectors that are horizontal with respect
to the connection and $\sigma(V)=1$.

We recall the definition of these elements. Take an open covering
$M=\cup U_\alpha$ where $K$ has the local unitary bases $b_\alpha$
and the transition functions
$$b_\beta=\gamma_{\alpha\beta}b_\alpha,\;\;\gamma_{\alpha\beta}=
e^{2\pi i\vartheta_{\alpha\beta}},$$ where we use $S^1=\{e^{2\pi
it}\,/\,t\in\mathds{R}\}$.

Then $\nabla$ has the local equations $\nabla b_\alpha=\omega_\alpha
b_\alpha$ where $\omega_\alpha$ are the local connection forms and
$\omega_\beta=\omega_\alpha+d\gamma_{\alpha\beta}.$ The preservation
of the Hermitian norm by $\nabla$ implies that $\omega_\alpha$ are
purely imaginary forms and we denote $\omega_\alpha=2\pi
i\varpi_\alpha$.

The above description corresponds to the Lie algebra identification
$ \mathit{u}(1)=span_{\mathds{R}}\{2\pi i\}$. If, instead, we take $
\mathit{u}(1)=\mathds{R}$, connection theory (e.g., see \cite{KN})
tells us that the form $\sigma$ is defined by the formula
\begin{equation}\label{exprsigma} \sigma|_{p^{-1}(U_\alpha)}
=p^*\varpi_\alpha+dt.\end{equation} Then, the curvature form of
the principal connection produced by $\nabla$ on $Q$ is
$\Omega=p^*(d\varpi_\alpha)$ and one has
\begin{equation}\label{Rprinsigma}
R_\nabla(X,Y)=-2\pi i\sigma([X^H,Y^H]),\end{equation} where the
upper index $H$ denotes the horizontal lift with respect to the
principal bundle connection. The horizontal lift $X^H(q)$, $q\in
Q$, is defined for all $X\in T_xM$ $(x\in M,\,p(q)=x)$ and it is
characterized by $p_*(X^H)=X,\;\sigma(X^H)=0.$ Finally, the
expression (\ref{exprsigma}) allows us to check that
$V={\partial}/{\partial t}$.

The importance of the prequantization space in geometric
quantization comes from the following result (see \cite{K} for the
symplectic case):
\begin{prop}\label{propK} Let $C^\infty_{inv}(Q,\mathds{C})$ be the space
of right-translation invariant, complex functions on $Q$. There
exists a natural isomorphism of complex linear spaces $\Gamma
K\approx C^\infty_{inv}(Q,\mathds{C})$ that transposes the action
of the quantum operator $\hat h$ $(h\in C^\infty_{wHam})$ to the
derivative defined by the vector field
\begin{equation}\label{hambar}
\bar X_h=X^H_h-[(\theta(\mathcal{X}_h)+h)\circ
p]V\in\chi^1(Q).\end{equation}\end{prop}
\begin{proof} The formula $s(p(q))=\bar s(q)q$, where
$q\in Q$ is seen as a basis of the fiber $K_{p(q)}$ and
$s\in\Gamma K$, $\bar s\in C^\infty_{inv}(Q,\mathds{C})$, defines
an isomorphism $\Gamma K\approx C^\infty_{inv}(Q,\mathds{C})$.
This isomorphism sends the function $X^H_h\bar s$ to the cross
section $\nabla_{X_h}s$ (Proposition III.1.3 of \cite{KN}) and a
calculation via coordinates $(x^i,t)$, where $(x^i)$ are
coordinates on $M$, shows that $V\bar s$ corresponds to $-2\pi
is$. Hence, the operator (\ref{KS}) is transformed into
(\ref{hambar}).\end{proof}

Proposition
\ref{propK} and property $\widehat{\{f,h\}}=[\hat f,\hat h]$ imply:
\begin{equation}\label{comutantLie} [\bar X_f,\bar X_h]=\bar
X_{\{f,h\}},\;\;\;\forall f\in C^\infty_{Ham}(M,E),h\in
C^\infty_{wHam}(M,E).\end{equation} For instance, for the constant
functions $f=0,f=1$ we may use $\mathcal{X}_f=(0,0)$, which gives
$\bar X_0=0,\bar X_1=-V$ and
\begin{equation}\label{Vauto} [X^H,V]\stackrel{(\ref{hambar})}{=}
-[V,\bar X_h]\stackrel{(\ref{comutantLie})}{=}\bar
X_{\{1,h\}}=\bar X_0=0.
\end{equation}

Like in \cite{WZ}, we produce a geometric structure that
accommodates the geometric elements defined by a g.p. data system
on the prequantization space $Q$. This structure is defined on the
stable tangent bundle \cite{V-stable}
$${\bf T}^{big}Q=T^{big}Q\oplus\mathds{R}^2=(TM\times\mathds{R})
\oplus(T^*M\times\mathds{R}).$$

For the vectors of ${\bf T}^{big}M$ we will use the following
notation (pay attention to boldface characters)
\begin{equation}\label{elementstabil}{\bf
X}=(\{X,u\},\{\alpha,v\})=(\mathcal{X},u,v),\;\;\mathcal{X}=(X,\alpha),
u,v\in\mathds{R}.\end{equation}
We recall that ${\bf T}^{big}M$ has the neutral metric
$${\bf{g}}({\bf{X}}_1,{\bf{X}}_2)
=\frac{1}{2}(\alpha_1(X_2)+\alpha_2(X_1)+u_1v_2+u_2v_1)$$ and the
{\it Wade bracket} \cite{Wd} \begin{equation}\label{Wade}
[(\{X_1,u_1\},\{\alpha_1,v_1\}),(\{X_2,u_2\},\{\alpha_2,v_2\})]_W\end{equation}
$$=(\{[X_1,X_2],X_1u_2-X_2u_1\},\{L_{X_1}\alpha_2-L_{X_2}\alpha_1+
\frac{1}{2}d(\alpha_1(X_2)-\alpha_2(X_1))$$
$$+u_1\alpha_2-u_2\alpha_1+\frac{1}{2}(v_2du_1-v_1du_2-u_1dv_2+u_2dv_1),$$
$$X_1v_2-X_2v_1+\frac{1}{2}(\alpha_1(X_2)-\alpha_2(X_1)-u_2v_1+u_1v_2)\}).$$

Formula (\ref{hambar}) suggests defining the {\it horizontal lift}
of the bundle $E$ by \cite{WZ}:
\begin{equation}\label{Ehorizontal}
E^H_q=\{(X^H-\theta(X,\alpha)V,p^*\alpha)\,/\,(X,\alpha)\in
E_{p(q)}\} \subseteq T_q^{big}Q\;(q\in Q).\end{equation} Obviously,
$E^H$ is a differentiable, big-isotropic bundle of the same rank as
$E$. The connection with (\ref{hambar}) is that every pair in
$E^H_q$ is the value of a pair $(\bar X_h,d(h\circ p))$ at $q$;
indeed, there exists a function $h$ such that
$(X_h,dh)_{p(q)}=(X,\alpha)$ and $h(p(q))=0$. Furthermore, we define
the {\it stable lift} ${\bf E}^{H}\subseteq{\bf T}^{big}Q$ as
follows
\begin{equation}\label{sEhorizontal}
{\bf E}^H_q=\{(\{X^H-\theta(X,\alpha)V,0\},
\{p^*\alpha,0\})\,/\,(X,\alpha)\in E_{p(q)}\}\oplus{\rm span}\{\bf V\}\end{equation}
where ${\bf V}=(\{V,0\},\{0,1\})$.
\begin{prop}\label{propWZ} ${\bf E}^H$ is an integrable, isotropic
subbundle of ${\bf T}^{big}Q$, ${\rm rank}{\bf E}^H={\rm rank}
E+1$, with the $\bf{g}$-orthogonal subbundle
\begin{equation}\label{sE'horizontal}
{\bf E}^{'H}_q=\{(\{Y^H-\theta(Y,\beta)V,0\},
\{p^*\beta,0\})\,/\,(Y,\beta)\in E'_{p(q)}\}\end{equation}
$$\oplus{\rm span}\{{\bf{V}}, (\{0,0\},\{0,1\}),{\bf
U}=(\{U^H,-1\},\{\sigma+p^*\nu,0\})\}.$$ \end{prop}
\begin{proof} The integrability of ${\bf E}^H$ means the closure of
$\Gamma{\bf{E}}^H$ with respect to the Wade bracket. We compute
the Courant bracket of two cross sections of $E^H$:
$$[(X_1^H-\theta(X_1,\alpha_1)V,p^*\alpha_1),
(X_2^H-\theta(X_2,\alpha_2)V,p^*\alpha_2)\stackrel{(\ref{Vauto})}{=}
([X_1^H,X_2^H]-X_1(\theta(X_2,\alpha_2))V$$
$$+X_2(\theta(X_1,\alpha_1))V, L_{X_1^H-\theta(X_1,\alpha_1)V}(p^*\alpha_2) -
L_{X_2^H-\theta(X_2,\alpha_2)V}(p^*\alpha_1)
+p^*(d(\alpha_1(X_2)))).$$ If we express the Lie derivative by the
Cartan formula $L_X=i(X)d+di(X)$, the $T^*M$-component becomes
$p^*(L_{X_1}\alpha_2-L_{X_2}\alpha_1+d(\alpha_1(X_2))$. Then,
since $[X_1,X_2]^H=pr_H[X_1^H,X_2^H]$, the $TM$-component is equal
to
$$[X_1,X_2]^H+\sigma([X_1,X_2])V-(d_{tr}\theta)((X_1,\alpha_1),(X_2,\alpha_2))V
-\theta([(X_1,\alpha_1),(X_2,\alpha_2)])V$$
$$\stackrel{(\ref{curburapct}),(\ref{Rprinsigma})}{=}
[X_1,X_2]^H-\theta([(X_1,\alpha_1),(X_2,\alpha_2)])V
+\omega_E((X_1,\alpha_1),(X_2,\alpha_2))V.$$

Now, for $\mathcal{X}_1,\mathcal{X}_2\in\Gamma E^H$ we get
$$\begin{array}{l}[(\mathcal{X}_1,0,0),(\mathcal{X}_2,0,0)]_W=
[\mathcal{X}_1,\mathcal{X}_2]_C+(\{0,0\},\{0,\alpha_1(X_2)\})\vspace*{2mm}\\
=(([X_1,X_2]^H-\theta([(X_1,\alpha_1),(X_2,\alpha_2)]_C)V,
p^*(L_{X_1}\alpha_2
-L_{X_2}\alpha_1\vspace*{2mm}\\+d(\alpha_1(X_2)))),0,0)+\alpha_1(X_2)(\{V,0\},\{0,1\})
\in\Gamma{\bf E}^H.\end{array}$$ The proof of the integrability of ${\bf E}$ is
completed by the simple calculation
$$[{\bf{V}},((X^H+\lambda V,p^*\alpha),0,0)]_W=\bf{0}.$$
The proof of the other assertions of the proposition is
straightforward.\end{proof}
\begin{rem}\label{obsWZ} {\rm Assume that the pair
$(U,\nu)$ is $g$-isotropic (e.g., $\nu=0$, see Proposition
\ref{corolarnu0}) and that $E$ is a Dirac structure.
Then, ${\bf{E}}^H\oplus{\rm span}\{{\bf U}\}$ is a {\it
Jacobi-Dirac structure} \cite{WZ}. Indeed, it is easy to get
$[{\bf{U}},{\bf{V}}]_W={\bf0}$. The only remaining condition
$[((X^H-\theta(X,\alpha)
V,p^*\alpha),0,0),{\bf{U}}]_W\in\Gamma{\bf{E}}^H$ can be deduced
from properties of the Wade bracket. The latter is conformally
related to a Courant bracket on $Q\times\mathds{R}$
\cite{V-stable} and, if one proceeds like in Remark 1.1 of
\cite{V-stable}, one gets
\begin{equation}\label{auxJD}{\bf g}([{\bf X},{\bf X}_1]_W,{\bf X}_2)+
{\bf g}({\bf X}_1,[{\bf X},{\bf X}_2]_W)=0,\;\forall{\bf X},{\bf
X}_1,{\bf X}_2\in\Gamma{\bf E}\end{equation} where ${\bf E}$ is
any almost Jacobi-Dirac structure on $Q$. Taking ${\bf
X}=((X^H-\theta(X,\alpha) V,p^*\alpha),0,0)$, ${\bf X}_1={\bf U}$
and, successively, ${\bf X}_2=((X^{'H}-\theta(X',\alpha')
V,p^*\alpha'),$ $0,0)$, ${\bf X}_2={\bf V}$, ${\bf X}_2={\bf U}$
in (\ref{auxJD}) we get $[((X^H-\theta(X,\alpha)
V,p^*\alpha),0,0),{\bf{U}}]_W\perp_{{\bf g}}{\bf E}^H$. Since we
are in the case where ${\bf E}^H$ is maximal ${\bf g}$-isotropic
we are done.}\end{rem}

Another interesting, but less comprehensive, structure on $Q$ is
defined by the $p$-pullback of the structure $E$ to $Q$, which
turns out to be:
\begin{equation}\label{pullback} (p^*E)_q= \{(X^H+\lambda
V,p^*\alpha)\,/\,(X,\alpha)\in E_{p(q)}\}=E^H\oplus
span\{(V,0)\}.\end{equation} Obviously, $p^*E$ is differentiable
and ${\rm rank}\,p^*E={\rm rank} E+1$. Integrability of $E$ and
(\ref{Vauto}) imply that the bracket of two cross sections of
$E^H$ belongs to $p^*E$. Since we also have
$$[(V,0),(X^H-\theta(X,\alpha)V,p^*\alpha)]=0,$$ we see that $p^*E$ is
closed by Courant brackets.

Since we have
$$(L_V(X^H-\theta(X,\alpha)V),L_V(p^*\alpha))=(0,0),\;(L_VV,L_V0)=(0,0),$$
the vector field $V$ is an infinitesimal automorphism of $p^*E$.
Accordingly, we may apply the prolongation construction of Theorem
2.1 of \cite{V-stable}, which gives the subbundle
\begin{equation}\label{eqprelung}
\tilde{E}=span\{(\{X^H+\lambda
V,0\},\{p^*\alpha,0\}),(\{0,0\},\{0,1\})\} \subseteq{\bf
T}^{big}M.\end{equation} A straightforward calculation shows that
$\tilde{E}$ is integrable too.
\begin{example}\label{bfPgraph} {\rm In the case of the structure
$E_{(\lambda,S)}$ of Example \ref{exgraphlambda} the pullback to
$Q$ is
$$p^*E_{(\lambda,S)}=graph(\flat_{p^*\lambda}|_{{\rm
span}\{Z^H,V\},Z\in S}).$$ In the case of the structure
$E_{(P,\Sigma)}$ of Example
\ref{exgraphP}, if we consider the bivector field $\Pi=P^H+V\wedge
U^H\in\chi^2(Q)$ and the morphism
$\Psi:T^*Q\times\mathds{R}\rightarrow TQ\times\mathds{R}$ defined
by
$$\Psi(\kappa,v)=(\sharp_\Pi\kappa+(v+P(\nu,\kappa))V,-\kappa(V)),\;\;
(\kappa\in T^*(Q),v\in\mathds{R}),$$ one has
$${\bf{E}}_{(P,\Sigma)}^H=graph(\Psi|_{p^*\Sigma\oplus{{\rm
span}}\{0,1\}}).$$}\end{example}
\section{Polarizations} In this section we discuss a problem that arises
in the comparison of geometric prequantization with quantization
commonly used in physics.

This discussion is motivated on one hand by the necessity to remove
the ambiguity of the quantum operator due to the non-uniqueness of
the (weak)-Hamiltonian vector field and on the other hand by the
following example.

The dynamics of a mechanical system with holonomic constraints may
be defined by a weak-Hamiltonian vector field with respect to an
integrable, big-isotropic structure of the type $E_{(P,\Sigma)}$.
Namely \cite{V-JMP}, assume that the configuration space of the
system is the manifold $N$ with local coordinates $(q^i)$, the
phase space is $M=T^*N$ with canonical, local coordinates
$(q^i,p_i)$ and the constraints are defined by the regular,
integrable distribution $L\subseteq TN$. In what follows we use
the Einstein summation convention.

Take $$P=\frac{\partial}{\partial
q^i}\wedge\frac{\partial}{\partial p_i},\;\Sigma={\rm
ann}\,\sharp_P(\pi^*{\rm ann}\,L),$$ where $\pi:M\rightarrow N$ is
the natural projection. Since $P$ is defined by the exact
symplectic form $\omega=dq^i\wedge dp_i\in\Omega^2(M)$, the
prequantization condition (\ref{integralptPSigma}) is satisfied.
More exactly, we may use a g.p. data system given by the trivial
bundle $K$ with basis $1$ and metric $\|1\|=1$, the connection
$\nabla$ defined by $\nabla1=2\pi i(p_idq^i)1$ and the cochain
$\theta=0$. Then, formula (\ref{eqE'P}) shows that the Hamiltonian
vector fields of the function $q^i$ are
$$X_{q^i}=\frac{\partial}{\partial p_i}+
\alpha_a\varphi^a_i\frac{\partial}{\partial p_i}$$ where
$\varphi^a_idq^i=0$ are the (independent) equations of $L$. The
corresponding quantum operator is
$$\hat q^i(s1)=[\frac{\partial s}{\partial
p_i}+\alpha_a\varphi^a_i\frac{\partial s}{\partial p_i}+ 2\pi
iq^is]1.$$ The result is unambiguous and reduces to multiplication
by $2\pi iq^i$ as required by physics\footnote{In fact, physics
requires multiplication by $q^i$, which happens if we divide the
Kostant-Souriau formula by $2\pi i$.} if $\partial s/\partial
p_i=0$.

In symplectic geometry, the distribution ${\rm
span}\,\{\partial/\partial p_i\}$ is called a polarization, and we
want a corresponding notion in the general case. We will extend
the definition that we gave in \cite{V-carte} for Poisson
manifolds.
\begin{defin}\label{defpol} {\rm A {\it real polarization} of an
integrable, big-isotropic structure $E$ is a pair of subspaces
$\mathcal{P}\subseteq\Gamma E,\mathcal{P}'\subseteq\Gamma E'$ with
the following properties
\begin{equation}\label{condpol}\begin{array}{c}
\mathcal{P}\subseteq\mathcal{P}',\;\chi^1(M)\cap\Gamma E\subseteq\mathcal{P},\;
\chi^1(M)\cap\Gamma E'\subseteq\mathcal{P}',\vspace*{2mm}\\

[\mathcal{P},\mathcal{P}]_C
\subseteq\mathcal{P},\;[\mathcal{P},\mathcal{P}']_C
\subseteq\mathcal{P}',\;\omega|_{\mathcal{P}\times\mathcal{P}'}=0.
\end{array}\end{equation} A {\it complex polarization} of $(E,E')$
is defined in the same way but replacing $E,E',\chi^1(M)$ by their
complexifications
$E_c=E\otimes\mathds{C},E'_c=E'\otimes\mathds{C},\chi^1_c(M)=\Gamma
T_cM$, $T_cM=TM\otimes\mathds{C}$.}\end{defin}

For the simplicity of notation we refer to real polarizations but
the results hold for complex polarizations as well.

With a polarization, we can associate the subspaces
$$\Gamma_{\mathcal{P}}K=\{s\in\Gamma
K\;/\,\nabla_Ys=-2\pi i\theta(Y,\beta)s, \;\forall
(Y,\beta)\in\mathcal{P}\},$$ $$\Gamma_{\mathcal{P}'}K=\{s\in\Gamma
K\,/\,\nabla_Zs=-2\pi i\theta(Z,\zeta)s, \,\forall
(Z,\zeta)\in\mathcal{P}'\},$$ which satisfy the condition
$\Gamma_{\mathcal{P}'}K\subseteq\Gamma_{\mathcal{P}}K$, as well as
the space of {\it polarized Hamiltonians}
$$C^\infty_{Ham}(M,\mathcal{P},\mathcal{P}')=\{f\in C^\infty_{Ham}(M,E)\,/\,
[(X_f,df),(Z,\zeta)]\in\mathcal{P}',
\,\forall (Z,\zeta)\in\mathcal{P}'\}$$ and the space
of {\it polarized weak-Hamiltonians}
$$C^\infty_{wHam}(M,\mathcal{P},\mathcal{P}')=\{h\in C^\infty_{wHam}(M,E)\,/\,
[(Y,\beta),(X_h,dh)]\in\mathcal{P}', \;\forall
(Y,\beta)\in\mathcal{P}\}.$$ In the case of a complex polarization
the previous spaces will be assumed to consist of complex valued
functions.
\begin{example}\label{expolar} {\rm The pair
$\mathcal{P}=\chi^1(M)\cap\Gamma
E,\mathcal{P}'=\chi^1(M)\cap\Gamma E'$ is a polarization such that
the corresponding spaces
$\Gamma_{\mathcal{P}}K,\Gamma_{\mathcal{P}'}K$ are
$$\Gamma_{E}K=\{s\in\Gamma
K\;/\,\nabla_Ys=-2\pi i\theta(Y,0)s, \;\forall Y\in
\chi^1(M)\cap\Gamma E\},$$
$$\Gamma_{E'}K=\{s\in\Gamma
K\,/\,\nabla_Zs=-2\pi i\theta(Z,0)s, \,\forall Z\in
\chi^1(M)\cap\Gamma E'\},$$ respectively, and
$$C^\infty_{Ham}(M,\mathcal{P},\mathcal{P}')=
C^\infty_{Ham}(M,E),\;
C^\infty_{wHam}(M,\mathcal{P},\mathcal{P}')=
C^\infty_{wHam}(M,E).$$ The restrictions $\hat
h|_{\Gamma_{E'}K},\hat f|_{\Gamma_{E}K}$ where $h\in
C^\infty_{wHam}(M,E),f\in C^\infty_{Ham}(M,E)$ are independent of
the chosen Hamiltonian vector fields. This polarization does not
solve the difficulty indicated by the example of the constrained
mechanical systems and there is a need for bigger, preferably
maximal, polarizations.}\end{example}
\begin{example}\label{expolar2} {\rm Let $E={\rm graph}\,\sharp_{\Pi}$
be the Dirac structure associated with a Poisson bivector field
$\Pi$. Then there exists a bijective correspondence between the
complex polarizations $(\mathcal{P},\mathcal{P}'=\mathcal{P})$ and
the subalgebras $\mathcal{Q}$ of the Lie algebra
$(\Omega^1\otimes\mathds{C},\{.,.\}_\Pi)$ (such a subalgebra
defined the notion of a polarization in \cite{V-carte}). This
correspondence is given by $\mathcal{P}\mapsto\mathcal{Q}={\rm
pr}_{\Omega^1(M)}\mathcal{P}$. Furthermore, the space of polarized
Hamiltonians is given by \cite{V-carte}
$$C^\infty_{Ham}(M,\mathcal{P},\mathcal{P})=\{f\in C^\infty(M,\mathds{C})\,/\,
\{df,\alpha\}_\Pi\in\mathcal{Q},
\,\forall\alpha\in\mathcal{Q}\}.$$
Finally, if $\Pi$ is prequantizable by a Hermitian line bundle $K$
and a contravariant derivative $D$ \cite{V-carte}, we may take an
arbitrary Hermitian connection $\nabla$ on $K$ and define a
cochain by the formula
$$\theta(\sharp_\Pi\xi,\xi)s=\frac{1}{2\pi i}(D_\xi
s-\nabla_{\sharp_\Pi\xi}s),\;\;\;s\in\Gamma K.$$ Then,
$(K,\nabla,\theta)$ is a g.p. data system and one has
\cite{V-carte}
$$\Gamma_{\mathcal{P}}K=\{s\in\Gamma K\,/\,D_\xi
s=0,\,\forall\xi\in\mathcal{Q}\}.$$ For instance, if
$M=\mathds{R}^{2n+h}=\{(q^i,p_i,t^u)\}$ and
$$\Pi=\sum_{i=1}^n\varphi_i(t^u)\frac{\partial}{\partial q^i}
\wedge\frac{\partial}{\partial p_i}$$ then $\mathcal{Q}={\rm
span}\{\varphi_i(t^u)dq^i\}$ (no summation) produces a
polarization $\mathcal{P}$.}\end{example}

Now, we can prove the following proposition.
\begin{prop}\label{valori} For any functions $h\in
C^\infty_{wHam}(M,\mathcal{P},\mathcal{P}'),f\in
C^\infty_{Ham}(M,\mathcal{P},\mathcal{P}')$, the operators $\hat
h$ restrict to a unique operator $\hat h
:\Gamma_{\mathcal{P}'}K\rightarrow\Gamma_{\mathcal{P}}K$ and the
operators $\hat f$ restrict to a unique operator $\hat
f:\Gamma_{\mathcal{P}'}K\rightarrow\Gamma_{\mathcal{P}'}K.$
\end{prop}
\begin{proof} The operators $\hat h|_{\Gamma_{\mathcal{P}'}K}$,
$\hat f|_{\Gamma_{\mathcal{P}'}K}$ are independent of the choice
of the Hamiltonian vector fields required by (\ref{KS}) because of
the second and third conditions (\ref{condpol}) and since
$\Gamma_{\mathcal{P}'}K\subseteq\Gamma_{\mathcal{P}}K$. Consider a
field $(Y,\beta)\in\mathcal{P}$ and a cross section
$s\in\Gamma_{\mathcal{P}'}K$. Then,
$$\nabla_Y(\hat hs)=\nabla_Y\nabla_{X_h}s
+2\pi iY[(\theta(\mathcal{X}_h)+h]s +2\pi
i[\theta(\mathcal{X}_h)+h]\nabla_Ys$$
$$=R_\nabla(Y,X_h)s+\hat h(\nabla_{Y}s) +\nabla_{[Y,X_h]}s+
2\pi i[Y(\theta(\mathcal{X}_h))+Yh]s.$$ In the result we insert
$$\nabla_Ys=-2\pi i\theta(Y,\beta)s,\;
\nabla_{[Y,X_h]}s=-2\pi i\theta([(Y,\beta),(X_h,dh)])s$$ (the second
equality is implied by the definition of
$C^\infty_{wHam}(M,\mathcal{P},\mathcal{P}')$) and
$$R_\nabla(Y,X_h))\stackrel{(\ref{curburapct})}
{=}-2\pi i[d_{tr}\theta((Y,\beta), (X_h,dh))-\omega_E((Y,\beta),
(X_h,dh))],$$ where, in fact,
$$\omega_E((Y,\beta), (X_h,dh))=-(Yh),$$ because the arguments are
$g$-orthogonal ($(Y,\beta)\in\mathcal{P}\subseteq\Gamma E,
(X_h,dh)\in\Gamma E'$). Then, after reductions, we get the
required result:
$$\nabla_Y(\hat hs)=-2\pi i\theta(Y,\beta)(\hat hs).$$ The
same calculation for $f\in
C^\infty_{Ham}(M,\mathcal{P},\mathcal{P}'),s\in\Gamma_{\mathcal{P}'}K,(Y,\beta)\in
\mathcal{P}'$ yields the second conclusion.\end{proof}
\begin{rem}\label{obslagrangiana} {\rm
The condition $\omega|_{\mathcal{P}\times\mathcal{P}'}=0$ that
enters in (\ref{condpol}) played no role in the previous proof.
However, it must be imposed because it is a necessary condition
for $\Gamma_{\mathcal{P}}K\neq0,\Gamma_{\mathcal{P}'}K\neq0$. This
follows by using (\ref{curburapct}) for
$(Y,\beta)\in\mathcal{P},(Z,\zeta)\in\mathcal{P}',
s\in\Gamma_{\mathcal{P}'}K\subseteq\Gamma_{\mathcal{P}}K$.
However, this condition is not sufficient and
$\Gamma_{\mathcal{P}}\neq0,\Gamma_{\mathcal{P}'}\neq0$ have to be
assumed.}\end{rem}

If the polarization is of the form $\mathcal{P}=\Gamma
P,\mathcal{P}'=\Gamma P'$ where $P,P'$ are subbundles of $E,E'$,
respectively, we have
$$\Gamma_{\mathcal{P}}K=\Gamma_{P}K=\{s\in\Gamma
K\,/\,\nabla_{Y_{x}}s=-2\pi i\theta(Y_x,\beta_x)s, \;\forall
(Y_{x},\beta_x)\in P_x, x\in M\},$$
$$\Gamma_{\mathcal{P}'}K=\Gamma_{P'}K=\{s\in\Gamma
K\,/\,\nabla_{Z_{x}}s=-2\pi i\theta(Z_x,\zeta_x)s, \;\forall
(Z_{x},\zeta_x)\in P'_x, x\in M\},$$ since
$(\nabla_Ys)(x)=\nabla_{Y_{x}}s$, etc. and any $(Y_x,\beta_x)\in
P_x,(Z_x,\zeta_x)\in P'_x$ have global, differentiable extensions
in $\Gamma P,\Gamma P'$. This point-wise setting can be extended
as follows.

Recall that the distribution $\mathcal{E}={\rm pr}_{TM}E$ is a
generalized foliation that defines the characteristic leaves
$\mathcal{L}$ of $E$ \cite{V-iso}. Any cross sections
$\mathcal{X}\in\Gamma(E|_{\mathcal{L}}),
\mathcal{Y}\in\Gamma(E'|_{\mathcal{L}})$
have differentiable extensions say, $\tilde{\mathcal{X}}\in
\Gamma E,\tilde{\mathcal{Y}}\in\Gamma E'$ to $M$ (at least locally) and we can
define a bracket
\begin{equation}\label{leafcroset}
\lfloor\mathcal{X},\mathcal{Y}\rceil=[\tilde{\mathcal{X}},
\tilde{\mathcal{Y}}]|_{\mathcal{L}} \in\Gamma(E'|_{\mathcal{L}}),\end{equation}
which is independent on the choice of the extensions. Indeed (like
in the proof of Theorem 2.1 of \cite{IKV}), it suffices to show
that the bracket vanishes for
$\tilde{\mathcal{Y}}=\sum_i\tilde{\lambda}^i\tilde{\mathcal{B}}_i$
where $\tilde{\mathcal{B}}_i$ is a local basis of $E'$ and
$\tilde{\lambda^i}|_{\mathcal{L}}=0$ (we do not have to consider a
similar $\tilde{\mathcal{X}}$ because $E\subseteq E'$.) Since for
$g$-orthogonal arguments the Courant bracket behaves like a Lie
algebroid bracket, we have
$$[\tilde{\mathcal{X}},\sum_i\tilde{\lambda}^i\tilde{\mathcal{B}}_i]=
\sum_i\{\tilde{\lambda}^i[\tilde{\mathcal{X}},\tilde{\mathcal{B}}_i]
+(({\rm
pr}_{TM}\tilde{\mathcal{X}})\tilde{\lambda}^i)\tilde{\mathcal{B}}_i\},$$
which has the zero restriction to $\mathcal{L}$ because ${\rm
pr}_{TM}\mathcal{X}\in\mathcal{L}$.

Accordingly, if we denote by $\bar{\Gamma}E,\bar{\Gamma}E'$ the
spaces of (possibly non differentiable) cross sections of $E,E'$
that have differentiable restrictions to each leaf $\mathcal{L}$
then, we get a bracket
$$\lfloor\bar{\mathcal{X}},\bar{\mathcal{Y}}\rceil\in\bar\Gamma
E',\;\forall\bar{\mathcal{X}}\in\bar{\Gamma}E,\bar{\mathcal{Y}}\in\bar{\Gamma}E',$$
which belongs to $\bar\Gamma E$ if
$\bar{\mathcal{X}},\bar{\mathcal{Y}}\in\bar{\Gamma}E$.

By generalized subbundles $P\subseteq E,P'\subseteq E'$ we will
understand fields of subspaces of the fibers of $E,E'$ such that
the restrictions to each leaf $\mathcal{L}$ are differentiable
vector bundles along $\mathcal{L}$. For instance, $TM\cap E,TM\cap
E'$ are generalized subbundles of $E,E'$, respectively. We shall
use the notation $\bar\Gamma P,\bar\Gamma P'$ in the same sense as
for $E,E'$.
\begin{defin}\label{polpct} {\rm A {\it real point-wise polarization} of an
integrable, big-isotropic structure $E$ is a pair of generalized
subbundles $P\subseteq E,P'\subseteq E'$ with the following
properties
\begin{equation}\label{pctpol}\begin{array}{c}P\subseteq P',\,TM\cap
E\subseteq P,\,TM\cap E'\subseteq P',\vspace*{2mm}\\
\lfloor\bar{\Gamma}P,\bar{\Gamma}P\rceil
\subseteq\bar{\Gamma}P,\,\lfloor\bar{\Gamma}P,\bar{\Gamma}P'\rceil
\subseteq\bar{\Gamma}P',\,\omega|_{P\times P'}=0.\end{array}
\end{equation} A {\it complex point-wise polarization} is defined
in the same way using the complexified bundles
$E_c,E'_c,T_cM$.}\end{defin}

If $(P,P')$ is a point-wise polarization, the spaces
$\Gamma_PK,\Gamma_{P'}K$ may still be defined and we may also
define the following spaces of functions
$$C^\infty_{Ham}(M,P,P')=\{f\in C^\infty_{Ham}(M,E)\,/\,
\lfloor(X_f,df),(Z,\zeta)\rceil\in\bar{\Gamma}P',
\,\forall (Z,\zeta)\in\bar{\Gamma}P'\},$$
$$C^\infty_{wHam}(M,P,P')=\{h\in C^\infty_{wHam}(M,E)\,/\,
\lfloor(Y,\beta),(X_h,dh)\rceil\in\bar{\Gamma}P',
\,\forall (Y,\beta)\in\bar{\Gamma}P\}.$$

With this notation we get
\begin{prop}\label{cazdelocalizare} For any function $f\in
C^\infty_{Ham}(M,P,P')$, the operators $\hat f$ restrict to a well
defined operator $\hat f:\Gamma_{P'}K\rightarrow\Gamma_{P'}K.$
\end{prop} \begin{proof}  The uniqueness of $\hat
f|_{\Gamma_PK\supseteq\Gamma_{P'}K}$ follows from the second
condition (\ref{pctpol}). Take $(Z_x,\zeta_x)\in P'_x$ $(x\in M)$
and extend it to a differentiable cross section
$(Z,\zeta)\in\Gamma(P'|_{\mathcal{L}_x})$, where $\mathcal{L}_x$
is the characteristic leaf of $E$ through $x$. Let
$(\tilde{Z},\tilde\zeta)\in\Gamma E$ be a further extension of
$(Z,\zeta)$ to a neighborhood of $x$ in $M$. Then, for
$s\in\Gamma_{P'}K$ one has
\begin{equation}\label{auxYx}\begin{array}{c}
\nabla_{Z_x}(\hat fs)=\nabla_{Z_x}\nabla_{X_f}s
+2\pi iZ_x[(\theta(\mathcal{X}_f))+f]s\vspace{2mm}\\
+2\pi
i[\theta(\mathcal{X}_f)+f](x)\nabla_{Z_x}s.\end{array}\end{equation}
On the other hand, one has
$$R_\nabla(Z_x,X_f(x))s=\nabla_{Z_x}\nabla_{X_f}s-
\nabla_{X_f(x)}\nabla_{\tilde{Z}}s-\nabla_{[\tilde Z,X_f](x)}s,$$
where the following happen: 1) since $X_f(x)$ is tangent to
$\mathcal{L}_x$, $\nabla_{X_f(x)}\nabla_{\tilde{Z}}s$ depends only
on $\nabla_{\tilde{Z}}s|_{\mathcal{L}_x}$, which is equal to
$-2\pi i\theta(Z,\zeta)s|_{\mathcal{L}_x}$ by the definition of
$\Gamma_{P'}K$; 2) since $[(\tilde
Z,\tilde\zeta),(X_f,df)](x)=-\lfloor
(X_f,df)|_{\mathcal{L}_x},(Z,\zeta)\rceil(x)$, the definitions of
$C^\infty_{Ham}(M,P,P')$ and $\Gamma_{P'}K$ imply
$$\nabla_{[\tilde Z,X_f](x)}s=2\pi i\theta_x(\lfloor
(X_f,df)|_{\mathcal{L}_x},(Z,\zeta)\rceil(x))s(x).$$ If these
results are inserted in (\ref{auxYx}) the same reductions like in
the proof of Proposition \ref{valori} hold and one gets $\hat
fs\in\Gamma_{P'}K$.\end{proof}
\begin{rem}\label{obspctDirac} {\rm
In the case of a Dirac structure it is natural to consider only
polarizations with $\mathcal{P}'=\mathcal{P}$, $P'=P$,
respectively. Accordingly the definitions will not refer to
$\mathcal{P}',P'$ any more and Proposition \ref{cazdelocalizare}
extends Lemma 6.1 of \cite{WZ}, which was proven differently
there.}\end{rem}
\begin{example}\label{expol1} {\rm Consider the integrable,
big-isotropic structure $E_{(\lambda,S)}$ of Example
\ref{exgraphlambda}, where $\lambda$ is a closed $2$-form and $S$
is an involutive subbundle of $TM$. If $S=TM$, $E$ is the Dirac
structure defined by the presymplectic form $\lambda$, $TM\cap
E={\rm ker}\,\lambda$ and one has only one characteristic leaf
$\mathcal{L}=M$. In this case, the examination of conditions
(\ref{pctpol}) is easy and shows that a point-wise polarization
with $P'=P$ may be identified with a vector subbundle $L={\rm
pr}_{TM}P$ such that ${\rm ker}\,\lambda\subseteq L$,
$\lambda|_L=0$ and $L$ is involutive. In the symplectic case,
${\rm ker}\,\lambda=0$ and, if we ask maximality of $P$, $L$ is an
involutive, Lagrangian subbundle as required by the classical
definition of a polarization of a symplectic manifold. In both
cases one gets
\begin{equation}\label{Haminex1}C^\infty_{Ham}(M,P,P')=\{f\in
C^\infty_{Ham}(M,E)\,/\,[X_f,Y]\in\Gamma L\}.\end{equation} Then,
if $\lambda$ is an integral form and we take the g.p. data system
$(K,\nabla,0)$ (see Example
\ref{condqex1}), we have
\begin{equation}\label{PKinex1} \Gamma_PK=\{s\in\Gamma
K\,/\,\nabla_Ys=0,\,\forall Y\in L\}.\end{equation}

For $S\subset TM$, the characteristic leaves are the leaves of the
foliation $S$, which is regular, hence, it is natural to look for
point-wise polarizations defined by regular subbundles $P,P'$.
Necessarily, $P={\rm graph}(\flat_\lambda|_L)$ where $L={\rm
pr}_{TM}P$ is an involutive subbundle of $S$ such that
$\lambda|_L=0$ and $TM\cap E_{(\lambda,S)}={\rm ker}\,\lambda\cap
S\subseteq L$. Formula (\ref{perpElambda}) implies $TM\cap
E'_{(\lambda,S)}=S^{\perp_\lambda}$ and $P'$ has to be a convenient
enlargement of $P$ such that $S^{\perp_\lambda}\subseteq P'$.

Assume that we are in the particular case where $\lambda|_S$ is
non degenerate, hence, it induces symplectic forms of the leaves
of $S$. Then, $S\cap S^{\perp_\lambda}=0$ and, since ${\rm
ker}\,\lambda\subseteq S^{\perp_\lambda}$, $TM\cap
E_{(\lambda,S)}=0$. If we start with a Lagrangian subfoliation $L$
of $(S,\lambda|_S)$, we obtain a subbundle $P={\rm
graph}(\flat_\lambda|_L)$ as required and the addition of
$$P'=\{(X+Y,\flat_\lambda X+\flat_\lambda Y+\gamma)\,/\,X\in
L,Y\in S^{\perp_\lambda},\gamma\in{\rm ann}\,S\}$$ (we might have
omitted $\flat_\lambda Y$ that belongs to ${\rm ann}\,S$, but, it
is convenient to keep it for the calculation that follows) gives a
polarization. Indeed, we have $S^{\perp_\lambda}\subseteq P'$
(take $X=0,\gamma=-\flat_\lambda Y$), $[\Gamma P,\Gamma
P]\subseteq\Gamma P$ and $\omega|_{P\times P'}=0$. For the closure
condition $[\Gamma P,\Gamma P']\subseteq\Gamma P'$, it suffices to
check
\begin{equation}\label{eq1inexpol1}\begin{array}{l}
[(X,\flat_\lambda X),(Y,\flat_\lambda
Y)]=([X,Y],\flat_\lambda[X,Y])\in\Gamma P',\vspace{2mm}\\

[(X,\flat_\lambda X),(0,\gamma)]=(0,L_X\gamma)\in\Gamma
P',\end{array}\end{equation} where $Y$ is $L$-foliated (since the
pairs in the left hand sides of (\ref{eq1inexpol1}) locally span
$\Gamma P'$). Then $[X,Y]\in L$ and the first relation
(\ref{eq1inexpol1}) holds with $P'$ replaced by $P$. Finally, it is
easy to check that $(0,L_X\gamma)\in {\rm ann}\,S\subseteq P'$,
which proves the second relation (\ref{eq1inexpol1}).

Formula (\ref{PKinex1}) obviously remains valid. Formula
(\ref{Haminex1}) remains valid too. Indeed, using the fact that the
infinitesimal transformation $X_f$ preserves $S$ we get
$$C^\infty_{Ham}(M,P,P')=\{f\in
C^\infty_{Ham}(M,E)\,/\,[X_f,Z+Y]\in\Gamma(L\oplus
S^{\perp_\lambda})\}$$ for $Z\in\Gamma L,Y\in\Gamma
S^{\perp_\lambda}$. Since $L_{X_f}\lambda=di(X_f)\lambda=d^2f=0$,
$X_f$ preserves $S^{\perp_\lambda}$ too and $[X_f,Y]\in\Gamma
S^{\perp_\lambda}$. Therefore, we remain with the condition
$[X_f,Z]\in\Gamma L$ that appears in (\ref{Haminex1}). Furthermore,
it follows straightforwardly that
$$C^\infty_{wHam}(M,P,P')=\{h\in
C^\infty_{wHam}(M,E)\,/\,[X,X_h]\in\Gamma(L\oplus
S^{\perp_\lambda}),\,\forall X\in\Gamma L\}.$$

The simplest example of the situation discussed above is given by
$$M=\mathds{R}^{2n+h}=\{(q^i,p_i,t^u)\},\;
\lambda=\sum_{i=1}^ndq^i\wedge dp_i,$$ $$S={\rm span}
\{\frac{\partial}{\partial q^i},\frac{\partial}{\partial p_i}\},
\;S^{\perp_\lambda}={\rm span}
\{\frac{\partial}{\partial t^u}\},\;
L={\rm span}\{\frac{\partial}{\partial p_i}\}.$$ Then,
$$P={\rm span}\{(\frac{\partial}{\partial p_i},-dq^i)\},\,
P'=P\oplus{\rm span}
\{(\frac{\partial}{\partial t^u},0)\}\oplus{\rm span}\{(0,dt^u)\}$$
define a polarization.}\end{example}
\begin{example}\label{polCRFK} {\rm This example extends the
classical notion of a K\"ahler polarization to generalized
geometry. We briefly recall the framework following
\cite{{Galt},{V-CRF}}. A classical metric F-structure on a manifold $M$
is a pair $(F,\gamma)$ where $\gamma$ is a Riemannian metric,
$F\in End\,TM$, $F^3+F=0$ and \begin{equation}\label{condFmetric}
\gamma(FX,Y)+\gamma(X,FY)=0\;\Leftrightarrow\;\flat_\gamma\circ
F=-^t\hspace{-1pt}F\circ\flat_\gamma\end{equation} ($t$ denotes
transposition). A generalized metric F-structure is an analogous
structure on $T^{big}M$ and is equivalent with a system that
consists of two classical metric F-structures with the same
metric, $(\gamma,F_+,F_-)$ and a $2$-form $\psi\in\Omega^2(M)$.
Then, there are two injections $j_\pm:TM\rightarrow T^{big}M$
defined by
\begin{equation}\label{defj} j_+(X)=(X,\flat_{\psi+\gamma}X),\;
j_-(X)=(X,\flat_{\psi-\gamma}X) \hspace{2mm}(X\in
TM)\end{equation} such that $ T^{big}M={\rm im}\,j_+\oplus{\rm
im}\,j_-$ , where ${\rm im}\,j_+\perp_g{\rm im}\,j_-$ and
\begin{equation}\label{ggamma}g((X,\flat_{\psi\pm\gamma}X),
(Y,\flat_{\psi\pm\gamma}Y))=\pm\gamma(X,Y).\end{equation} The
injections $j_\pm$ send $F_\pm$ to structures $\mathcal{F}_\pm$ on
their images with $\mathcal{F}=\mathcal{F}_++\mathcal{F}_-\in
End(T^{big}M)$ such that $\mathcal{F}^3+\mathcal{F}=0$.

The structures $F_\pm$ yield decompositions
$$T_cM=H_\pm\oplus\bar H_\pm\oplus Q_{\pm c}\hspace
{3mm}(Q_{\pm c}=Q_\pm\times\mathds{C},\,Q_\pm\subseteq TM)$$ where
the terms are the $(\pm i,0)$-eigenbundles, respectively. If we
denote $E_\pm=j_\pm(H_\pm),S_\pm=j_\pm(Q_\pm)$,
(\ref{condFmetric}) and (\ref{ggamma}) imply that
\begin{equation}\label{defEinCRFK} E_1=E_+\oplus E_-,\;E_2=E_+\oplus\bar{E}_-
\end{equation} define complex, generalized, big-isotropic
structures with the $g$-orthogonal bundles
\begin{equation}\label{defE'inCRFK} E'_1=E_1\oplus
S_c,\;\;E'_2=E_2\oplus S_c\;(S=S_+\oplus
S_-,\,S_c=S\otimes\mathds{C}).\end{equation} The generalized
metric F-structure defined by $(\gamma,F_\pm,\psi)$ is said to be
an integrable or a CRFK-structure if the following Courant bracket
closure conditions are satisfied
\begin{equation}\label{CRFK1'} \begin{array}{l}[E_+,E_+]
\subseteq E_+,\,[E_+,S_{+c}]\subseteq E_+ \oplus S_{+c},\vspace{2mm}\\

[E_-,E_-]\subseteq E_-,\, [E_-,S_{-c}]\subseteq E_- \oplus S_{-c}.
\end{array}\end{equation} Equivalently, a CRFK-structure is
characterized by the fact that the structures $E_1,E_2$ are
integrable and $[S_+,S_-]\subseteq S$ \cite{V-CRF}. If $F_\pm$ are
complex structures $S=0$ and the CRFK-structure is a generalized
K\"ahler structure \cite{Galt}.

From (\ref{defEinCRFK}), (\ref{defE'inCRFK}) and (\ref{CRFK1'}),
we see that $(P=E_+,P'=E_+\oplus S_{\pm c})$ might define a
polarization of both $E_1$ and $E_2$. But, the algebraic
conditions included in (\ref{pctpol}) may not hold. Using
(\ref{defj}), it follows that
\begin{equation}\label{condpsi1} \psi(F_\pm X,Y)+\psi(X,F_\pm
Y)=0\end{equation} implies $\omega|_{P\times P'}=0$. Furthermore, it
is a technical matter to check that, if the tensors $\psi\pm\gamma$
are non degenerate (e.g., if $\psi=t\psi'$ with
$t\in\mathds{R}_{\geq0}$ small), then $TM\cap E_{1,2}=TM\cap
E'_{1,2}=0$. Hence, if $\psi$ is such that $\psi\pm\gamma$ are non
degenerate and (\ref{condpsi1}) holds the subbundles $P,P'$ define
polarizations of $E_1,E_2$. For a classical K\"ahler structure $E_1$
is determined by the antiholomorphic tangent bundle and $E_2$ is
defined by the K\"ahler form. Classically, the antiholomorphic
tangent bundle is regarded as a polarization of the symplectic
structure defined by the K\"ahler form.}\end{example}
\section{Appendix: Truncated cohomology}
Truncated cohomology may also be useful in other situations,
particularly in foliation theory.
\begin{defin}\label{struncform} {\rm Let $F\subseteq L$ be a
pair of Lie algebroids  over a manifold $M$. An $s$-{\it truncated
form of degree $k$} $(0\leq s\leq{\rm dim}\,L)$ is a multilinear
bundle morphism
\begin{equation}\label{s-forme} \lambda:\underbrace{L\times\cdots\times L}_{s-{\rm
times}}\times \underbrace{F\times\cdots\times F}_{(k-s)-{\rm times}}
\rightarrow(M\times\mathds{R})\end{equation} that is skew-symmetric
with respect to all the arguments in $F$ (even if their number is
larger than $k-s$) and with respect to all the arguments in
$L\backslash F$ separately.}\end{defin}

We will denote by $\Omega^k_s(L,F)$ the space of $s$-truncated
forms of degree $k$; if $k>s+{\rm dim}\,F$ then
$\Omega^k_s(L,F)=0$. The truncated forms may be seen as
restrictions of forms
$\tilde\lambda\in\Omega^k(L)=\Gamma\wedge^kL^*$. Indeed, choose a
complementary subbundle $Q$ of $F$ in $L$ $(L=Q\oplus F)$ and
consider the decomposition
\begin{equation}\label{descomp-p,q}
\wedge^kL^{*}=\oplus_{p+q=k}(\wedge^pQ^*)\otimes(\wedge^qF^{*})\;
(k=1,...,{\rm rank}\,L).\end{equation} Then,
$\forall\lambda\in\Omega^k_s(L,F)$ we get uniquely defined forms
$\lambda_{(p,q)}$ of {\it type} $(p,q)$ $(p+q=k)$ with $p\leq s$ by
evaluating $\lambda$ on $p$ arguments in $Q$ and $q$ arguments in
$F$ and if we take
\begin{equation}\label{formaindices} \tilde{\lambda}=\lambda_{(0,k)}
+...+\lambda_{(s,k-s)}\end{equation} $\lambda$ is the restriction
of $\tilde{\lambda}$. Hence, if we denote
$\Omega^{(p,q)}(Q,F)=\Gamma(\wedge^pQ^*\otimes\wedge^qF^{*})$, we
may identify
\begin{equation}\label{struncate} \Omega^k_s(L,F)\approx\tilde{\Omega}^k_s(Q,F)
=\oplus_{p+q=k,p\leq s}\Omega^{(p,q)}(Q,F).\end{equation}

We define the exterior differential $d^s_L$ of an $s$-truncated
form by the same formula as the one used for the usual exterior
differential $d_L$ with the last $\geq k-s+1$ arguments in $\Gamma
F$. If $\lambda\in\Omega^k_s(L,F)$ then
$d^s_L\lambda\in\Omega^{k+1}_s(L,F)$ and $(d^s_L)^2\lambda=0$ is a
consequence of the Jacobi identity for $L$.

Thus, $(\Omega_s^k(L,F),d^s_L)$ is a cochain complex and defines
cohomology spaces $H^k_s(L,F)$ that we call $s$-{\it intermediate
de Rham cohomology spaces}. For $s=0$ the result is the de Rham
cohomology of $F$ and for $s={\rm dim}\,L$ the result is the de
Rham cohomology of $L$. Restriction to arguments as required
yields homomorphisms
\begin{equation}\label{cohmorf} h_{s>u}:H^k_s(L,F)\rightarrow
H^k_u(L,F).\end{equation}

For the interpretation of the intermediate cohomology via
(\ref{struncate}) notice the following decomposition of the exterior
differential $d_L$:
\begin{equation}\label{descompdL}
d_L=(d'_L)_{(1,0)} + (d''_L)_{(0,1)}
+(\partial_L)_{(2,-1)},\end{equation} where the indices denote the
$(Q,F)$-type of the operators and the following relations hold
\begin{equation}\label{propdL} \begin{array}{c}(d''_L)^2=0, d'_Ld''_L +
d'_Ld''_L=0,\;(\partial_L)^2=0,\vspace{2mm}\\ d'_L\partial_L +
\partial_Ld'_L=0,\;(d'_L)^2 + d''_L\partial_L +
\partial_Ld''_L=0 \end{array}\end{equation} (see \cite{V73} for the
case where $L=TM$ and $F$ is a foliation on $M$).

Furthermore, put $\tilde{d}^s_L={\rm
pr}_{\tilde{\Omega}_s^{k+1}(Q,F)}\circ d_L$, which means that for
$\tilde{\lambda}$ defined by (\ref{formaindices}) one has
\begin{equation}\label{exprds} \tilde{d}^s_L\tilde{\lambda} = d_L\tilde\lambda
-d'_L\lambda_{(s,k-s)}-\partial_L\lambda_{(s,k-s)}
-\partial_L\lambda_{(s-1,k-s+1)}.\end{equation} If we use
(\ref{exprds}) in order to compute $(\tilde{d}^s_L)^2$,
(\ref{propdL}) shows that everything cancels. Accordingly,
$(\tilde{\Omega}^k_s(Q,F),\tilde{d}^s_L)$ is a cochain complex,
which is isomorphic with the $s$-{\it intermediate de Rham
complex} of the pair $(L,F)$ and produces isomorphic cohomology
spaces $H^k_s(L,F)$.

As an example, we prove a result for the intermediate cohomology
of a foliation $F\subseteq L=TM$; in this case we shall omit the
index $L$ in the notation of the various differentials involved.
\begin{prop}\label{Poincare}
Let $\lambda\in\Omega^k_s(TM,F)$ be defined by
(\ref{formaindices}) and assume that $d^s\lambda=0$. Then,
locally, if $k>s$ one has $\lambda=d^s\mu$ and if $k=s$ one has
$\tilde{\lambda}=d\tilde{\mu}+\tilde{\nu}$ where $\tilde{\nu}\in
\Gamma\wedge^sQ^*$ and $d''\tilde{\nu}=0$.\end{prop}
\begin{proof} In the proposition tilde has the significance given by
formula (\ref{formaindices}). However, for simplicity we omit the
sign tilde hereafter. We have a decomposition
\begin{equation}\label{auxappendix}
d^s\lambda=\sum_{j=0}^s\xi_{(j,k+1-j)}\end{equation} where
\begin{equation}\label{auxappendix1}
\xi_{(j,k+1-j)}=d''\lambda_{(j,k-j)}+d'\lambda_{(j-1,k+1-j)}
+\partial\lambda_{(j-2,k+2-j)}\end{equation} and $d^s\lambda=0$ is
equivalent to $\xi_{(j,k+1-j)}=0$ for $j=0,...,s$.

If $k>s$ we can prove by induction that there are local forms
$\mu$ such that
\begin{equation}\label{auxappendix2}
\lambda_{(j,k-j)}=d''\mu_{(j,k-1-j)} + d'\mu_{(j-1,k-j)}
+\partial\mu_{(j-2,k+1-j)},\;\;(j=0,...,s).\end{equation} For
$j=0$, (\ref{auxappendix2}) holds since
$$\xi_{(0,k+1)}=d''\lambda_{(0,k)}=0\;\Rightarrow\;
\lambda_{(0,k)}=d''\mu_{(0,k-1)},$$ for a local form $\mu$, in
view of the Poincar\'e lemma for the operator $d''$ (e.g.,
\cite{V73}).

Then, if we assume that (\ref{auxappendix2}) holds for lower
values of $j$, the annulation of the form given by
(\ref{auxappendix1}) together with the properties (\ref{propdL})
give
$$d''(\lambda_{(j,k-j)} + d'\mu_{(j-1,k-j)}
+\partial\mu_{(j-2,k+1-j)})=0$$ and the $d''$-Poincar\'e lemma
produces the local form $\mu_{(j,k-1-j)}$ such that
(\ref{auxappendix2}) holds for the index $j$.

Formulas (\ref{auxappendix}), (\ref{auxappendix1}) applied to the
form $\mu=\sum_{j=0}^s\mu_{(j,k-1-j)}$ yield the required
conclusion for $k>s$.

If $k=s$, formula (\ref{auxappendix2}) similarly holds up to
$j=s-1$, which implies
$$\lambda=\lambda_{(s,0)}+d\mu-d'\mu_{(s-1,0)}-\partial\mu_{(s-2,1)}$$
for a local form $\mu=\sum_{h=1}^{s-1}\mu_{(h,s-h-1)}$. On the
other hand, if we write down (\ref{auxappendix1}) for $j=s=k$
while inserting the values of
$\lambda_{(s-1,1)},\lambda_{(s-2,2)}$ given by
(\ref{auxappendix2}) and use (\ref{propdL}), we get
$$d''(\lambda_{(s,0)}-d'\mu_{(s-1,0)}-\partial\mu_{(s-2,1)})=0.$$
Therefore, if we denote
$$\nu=\lambda_{(s,0)}-d'\mu_{(s-1,0)}-\partial\mu_{(s-2,1)}
\in\Gamma\wedge^{s}Q^*,$$ we exactly have the required conclusion.
\end{proof}

Let $\Phi_s(F)$ be the sheaf of germs of forms of the type
$d\mu+\nu$, $\mu\in\Omega^{s-1}(M),\nu\in\wedge^{s}Q^*,d''\nu=0$
(i.e., $\nu$ is a basic $s$-form). Then, we get
\begin{corol}\label{corolA} The cohomology of $M$ with values in
$\Phi_s(F)$ is given by the formula
$$H^u(M,\Phi_s(F))=H^{s+u}_s(TM,F).$$\end{corol}
\begin{proof} Proposition \ref{Poincare} shows that the sequence
of sheaves
$$0\rightarrow \Phi^s(F)\stackrel{\subseteq}{\rightarrow}
\underline{\Omega}^s(M) \stackrel{d^s}{\rightarrow}\underline{\Omega}^{s+1}_s(TM,F)
\stackrel{d^s}{\rightarrow}\cdots$$
is a fine resolution of the sheaf $\Phi_s(F)$.\end{proof}

Proposition \ref{Poincare} and Corollary \ref{corolA} hold for the
more general case of pairs of Lie algebroids $F\subseteq L$ that
satisfy the relative Poincar\'e lemma defined below.
\begin{defin}\label{algLP} {\rm  1) A Lie algebroid $L$
{\it satisfies the Poincar\'e lemma} if, $\forall k\geq0$,
$\forall\lambda\in\Omega^{k+1}(L)$ with $d_L\lambda=0$, there
exist local forms $\mu\in\Omega^{k}(L)$ such that
$\lambda=d_L\mu$. 2) A pair of Lie algebroids $F\subseteq L$ {\it
satisfies the relative Poincar\'e lemma} if $\forall
p\geq0,q\geq0$, $\forall\lambda\in\Omega^{(p,q+1)}(Q,F)$ with
$d''_L\lambda=0$, there exist local forms
$\mu\in\Omega^{(p,q)}(L)$ such that
$\lambda=d''_L\mu$.}\end{defin}

If $F=L$ the relative Poincar\'e lemma is the Poincar\'e lemma for
$L$. The relative Poincar\'e lemma is independent of the choice of
$Q$ since a change of $Q$ only adds terms with more than $k-s$
arguments in $F$, which vanish. The operator $d''$ of a foliation
satisfies the relative Poincar\'e lemma but, unfortunately, we
know of no other interesting examples. (In the case of a complex
manifold the sheaf $\Phi_s$ coincides with the sheaf of
holomorphic $s$-forms and Proposition \ref{Poincare} gives nothing
new.)
\hspace*{7.5cm}{\small \begin{tabular}{l} Department of
Mathematics\\ University of Haifa, Israel\\ E-mail:
vaisman@math.haifa.ac.il \end{tabular}}

\begin{thebibliography}{xx}
\bibitem{C} T. J. Courant, Dirac Manifolds, Transactions Amer.
math. Soc., 319 (1990), 631-661.
\bibitem{Galt} M. Gualtieri, Generalized complex geometry, Ph.D.
thesis, Univ. Oxford, 2003; arXiv:math.DG/0401221.
\bibitem{IKV} V. Itskov, M. Karasev and Yu. M. Vorobjev,
Infinitesimal Poisson cohomology, A.M.S. Transl., 187 (2) (1998),
327-360.
\bibitem{KN} S. Kobayashi and K. Nomizu, Foundations of
Differential Geometry, Vol. I, Interscience Publ., New York, 1963.
\bibitem{K} B. Kostant, Quantization and unitary
representations, Lect. Notes in Math. 170, Springer-Verlag, New
York, 1970, 237-253.
\bibitem{Sour} J. M. Souriau, Structure des syst\`emes dynamiques,
Dunod, Paris, 1969.
\bibitem{V73} I. Vaisman, Cohomology and Differential Forms, M.
Dekker, Inc., New York 1973.
\bibitem{V-Torino} I. Vaisman, Basic ideas of geometric
quantization, Rend. Mat. Univ. Politecn. Torino, 37 (1979), 31-41.
\bibitem{V-carte} I. Vaisman, Lectures on the geometry of Poisson
manifolds, Progress in Math., vol. 118 Birkh\"auser Verlag, Boston,
1994.
\bibitem{V-stable} I. Vaisman, Dirac Structures and Generalized
Complex Structures on $TM\times\mathds{R}^h$, Advances in Geom. 7
(2007), 453-474.
\bibitem{V-iso} I. Vaisman, Isotropic Subbundles of $TM\oplus
T^*M$, Int. J. of Geom. Methods in Modern Phys., 4 (3) (2007),
487-516. \bibitem{V-JMP} I. Vaisman, Weak-Hamiltonian dynamical
systems, J. of Math. Phys., 48, 082903 (2007).
\bibitem{V-CRF} I. Vaisman, Generalized CRF-structures, Geometriae
Dedicata, http://dx.doi.org/10.1007/s10711-008-9239-z;
arXiv:0705.3934.
\bibitem{Wd} A. Wade, Conformal Dirac structures, Lett. Math. Phys.
53 (2000), 331-348.
\bibitem{WZ} A. Weinstein and M. Zambon, Variations on
Prequantization, Travaux mathematiques, Univ\'ersit\`e de
Luxembourg, Fascicule XVI (4$^{\rm th}$ Conference on Poisson
Geometry), pp. 187-219 (2005).
\end{thebibliography}
\end{document}